\documentclass[a4paper,12pt]{article}
\usepackage{times, url}
\usepackage{amsmath,amsthm} 
\usepackage[lite]{amsrefs}
\usepackage{chngcntr}
\usepackage[english]{babel}
\usepackage{enumitem}
\usepackage{hyperref} 
\usepackage{amsthm}
\usepackage{amssymb}
\usepackage[all,cmtip]{xy}
\newcommand{\legendre}[2]{\ensuremath{\left( \frac{#1}{#2} \right) }}

\usepackage{amsmath,amssymb,amsthm,amsfonts}

\newtheorem{theorem}{Theorem}[section]
\newtheorem{lemma}[theorem]{Lemma}

\newtheorem{corollary}[theorem]{Corollary}
\theoremstyle{definition}

\newtheorem{remark}{Remark}
\newtheorem*{case1*}{Case 1}
\newtheorem*{case2.a*}{Case 2.a}
\newtheorem*{case2.b*}{Case 2.b}
\newtheorem*{case2.c*}{Case 2.c}
\begin{document}

\setcounter{page}{1}

\begin{center}
{\LARGE \bf  Some new primality criteria based on Lucas sequences.}
\vspace{8mm}

{\large \bf Gaitanas Konstantinos}
\vspace{3mm}

Department of  Applied Mathematical and Physical Sciences\\
 National Technical University of Athens \\ 
Heroon Polytechneiou Str., Zografou Campus, 15780 Athens, Greece \\ 
e-mail: \url{kostasgaitanas@gmail.com}
\vspace{2mm}

\end{center}
\vspace{10mm}

\noindent
{\bf Abstract:} In this paper, we provide some novel results concerning the behavior of $\frac{U_{kn}}{U_k}$ modulo ${U_n}$, where $(U_n)_{n\in\mathbb{N}}$ is the Lucas sequence of the first kind. As a consequence, we obtain some primality criteria which do not seem to appear in the literature.
\vspace{10mm}
\section{Introduction} 
Let $a_n(P, Q)$ be numerical sequences, satisfying the recurrence relation 
\begin{gather*}a_{n+1}=P\cdot a_{n}-Q\cdot a_{n-1}\quad P, Q \in \mathbb{Z}.\end{gather*}The two special cases with initial conditions $(a_0, a_1)=(0, 1)$ and $(a_0, a_1)=(2, P)$ are the \emph{Lucas sequences} of the first and second kind and will be denoted by $U_n(P, Q)$ and $V_n(P, Q)$, respectively. Since the Lucas sequences satisfy a linear recurrence relation with constant coefficients, it is easy to obtain explicit expressions (Binet forms) from the characteristic polynomial $p(x)=x^2-Px+Q$. If the discriminant of $p(x)$ is non-zero and $\alpha, \beta$ are its roots, then 
\begin{gather*}U_n=\frac{\alpha^n-\beta^n}{\alpha-\beta}, \quad V_n=\alpha^n+\beta^n.\end{gather*} 
Many popular examples are recovered as particular instances of Lucas sequences, which are important for their own sake. Probably the two most popular are the Fibonacci and Lucas numbers (\emph{``The two shining stars in the vast array of integer sequences'' \cite{1}}) which are obtained from $U_n(1, -1)$ and $V_n(1, -1)$, respectively. Other well-known examples include the Pell, Mersenne, Jacobsthal and Balancing numbers. For more material that might satisfy the reader's curiosity we refer to the book of Paulo Ribenboim \emph{``My Numbers, My Friends''} \cite{2}. The main purpose of this paper is to present some congruence relations which hold among the terms of $(U_n)_{n\in\mathbb{N}}$ and eventually provide some primality criteria which are aimed to serve mostly for theoretical purposes. We emphasize that the congruences we obtain throughout this paper seem to be of independent interest \footnote{It is possible that some of these congruences have already appeared somewhere in the vast mathematical literature. For the convenience of the reader, we provide detailed proofs to make the paper self-contained.}
\section{Notation and Preliminaries}
This section is a quick-reference guide to the notation and background information that will be assumed throughout this paper.\\
We denote by $\varphi(n)$ Euler's totient function, which counts the positive integers not greater than $n$, that are relatively prime to $n$. We will also use without proof the fact that the number of positive integers $k$ less than $n$ with $\gcd(n, k)=d$, is exactly $\varphi(\frac{n}{d})$. Legendre's symbol is defined as \begin{gather*}\legendre{a}{p}=\begin{cases}
 1, & \text{if } a \text{ is a quadratic residue modulo } p \text{ and } a \not\equiv 0\pmod p, \\
-1, & \text{if } a \text{ is a quadratic nonresidue modulo } p, \\
 0, & \text{if } a \equiv 0 \pmod p.
\end{cases}.\end{gather*}In addition, we mention some well-known facts about Lucas sequences. For the sake of brevity, we omit their proofs since they are obtained without much effort and can be found in various mathematical contexts, such as \cite{2}. 
\begin{theorem}
\counterwithin{enumi}{theorem}
Let $D=P^2-4Q$. The following hold true:
\begin{enumerate}
\item\label{part1}$U_{m+n}=U_mV_{n}-Q^nU_{m-n}$.
\item\label{part2}$V_{n}=2U_{n+1}-PU_n$.
\item\label{part3}$V_{n}^2-DU_n^2=4Q^n$.
\end{enumerate}
\end{theorem}
Moreover, for the rest of the paper we will assume that $P$ and $Q$ are coprime. In this way, we ensure that $(U_n)$ is a \emph{strong divisibility sequence} which means that $\gcd(U_x, U_y)=U_{\gcd(x, y)}$. This evidently implies that if $y\mid x$, then $U_y\mid U_x$. This is something which will be frequently used. As for other prerequisites, the reader is expected to be familiar only with basic congruence rules and Bézout's identity, which states that if $\gcd(x, y)=d$, then there are $s, q\in \mathbb{Z}$ such that $sx+qy=d$. 
\section{Two important lemmas.}
In this section, we provide two lemmas which will be the key tools for the proof of Theorem \eqref{theor}. We note that a version of the first lemma already appears in \cite{3}, in the proof of Theorem 4.1.
\begin{lemma} Let $k\in \mathbb{N^*}$. The following holds true:
\begin{equation}\label{eq1}\frac{U_{kn}}{U_k}\equiv\begin{cases} 
      n\cdot Q^{\frac{k(n-1)}{2}}\pmod{U_k}, & \text{if}\quad n\equiv 1\pmod{2}, \\
      nU_{k+1}\cdot Q^{\frac{k(n-2)}{2}}\pmod{U_k}, & \text{if}\quad n\equiv 0\pmod{2}.\\
   \end{cases}\end{equation} 
\begin{proof}
We prove this by applying induction on $n$. We note that \eqref{eq1} obviously holds true for $n=0$ and $n=1$. Suppose that it holds true
for every natural number not greater than $n$, with $n\geq2$. From \eqref{part1} for $m=kn$ and $n=k$, we obtain 
\begin{equation}\label{eq2}\frac{U_{k(n+1)}}{U_k}=\frac{U_{kn}}{U_k}V_k-Q^k\frac{U_{k(n-1)}}{U_k}.\end{equation}
Assume that $n+1$ is odd. Then $n$ is even and $n-1$ is also odd. Hence, using the induction hypothesis the last equation yields 
\begin{gather*}\frac{U_{k(n+1)}}{U_k}\equiv nU_{k+1}Q^{\frac{k(n-2)}{2}}V_k-Q^k(n-1)Q^{\frac{k(n-2)}{2}}\pmod{U_k}.\end{gather*}
From \eqref{part2}, $V_k\equiv 2U_{k+1}\pmod{U_k}$. Thus, $nU_{k+1}V_k\equiv \frac{n}{2}\cdot 2U_{k+1}V_k\equiv\frac{n}{2}V_k^2\equiv 2nQ^k\pmod{U_k}$. The last congruence holds, since from \eqref{part3}, $V_k^2\equiv 4Q^k\pmod{U_k}$. It follows that 
\begin{gather*}\frac{U_{k(n+1)}}{U_k}\equiv 2nQ^kQ^{\frac{k(n-2)}{2}}-(n-1)Q^{\frac{k(n-2)}{2}}Q^k\equiv (n+1)Q^{\frac{kn}{2}}\pmod{U_k}\end{gather*} 
which proves the claim if $n+1$ is odd.\\
If $n+1$ is even, then $n$ is odd and $n-1$ is even. From the induction hypothesis, \eqref{eq2} yields 
\begin{gather*}\frac{U_{k(n+1)}}{U_k}\equiv n\cdot Q^{\frac{k(n-1)}{2}}V_k-Q^k(n-1)U_{k+1}Q^{\frac{k(n-3)}{2}}\pmod{U_k}.\end{gather*}
But $V_k\equiv 2U_{k+1}\pmod{U_k}$. Therefore, 
\begin{gather*}\frac{U_{k(n+1)}}{U_k}\equiv 2n\cdot Q^{\frac{k(n-1)}{2}}U_{k+1}-Q^k(n-1)U_{k+1}Q^{\frac{k(n-3)}{2}}\equiv(n+1)U_{k+1}Q^{\frac{k(n-1)}{2}}\pmod{U_k}\end{gather*} which proves the claim if $n+1$ is even.
\end{proof}
\end{lemma}
We use similar arguments to obtain our second result.
\begin{lemma}Let $k\in \mathbb{N}$. The following holds true:
\begin{equation}\label{eq3}U_{kn+1}\equiv\begin{cases} 
      U_{k+1}Q^{\frac{k(n-1)}{2}}\pmod{U_k}, & \text{if}\quad n\equiv 1\pmod{2}, \\
      Q^{\frac{kn}{2}}\pmod{U_k}, & \text{if}\quad n\equiv 0\pmod{2}.\\
   \end{cases}\end{equation} 
\begin{proof}
We prove this again by applying induction on $n$. It is evident that \eqref{eq3} holds true for $n=0$ and $n=1$. Suppose that it holds true
for every natural number not greater than $n$, with $n\geq2$. From \eqref{part1} for $m=kn$ and $n=k+1$, we obtain 
\begin{equation}\label{eq4}U_{k(n+1)+1}=U_{kn}V_{k+1}-Q^{k+1}U_{k(n-1)-1}\equiv-Q^{k+1}U_{k(n-1)-1}\mod{U_k}\end{equation} since $U_k\mid U_{kn}$. From definition, $U_{k(n-1)+1}=PU_{k(n-1)}-QU_{k(n-1)-1}$. But $U_k\mid U_{k(n-1)}$, which implies that $-QU_{k(n-1)-1}\equiv U_{k(n-1)+1}\pmod{U_k}$. Therefore, from \eqref{eq4} we deduce that 
\begin{equation}\label{eq5}U_{k(n+1)+1}\equiv Q^kU_{k(n-1)+1}\pmod{U_k}.\end{equation}
Assume that $n+1$ is even. Then $n-1$ is also even and from the induction hypothesis the last congruence yields $U_{k(n+1)+1}\equiv Q^kQ^{\frac{k(n-1)}{2}}\equiv Q^{\frac{k(n+1)}{2}}\pmod{U_k}$.\\
If $n+1$ is odd, then $n-1$ is also odd. From the induction hypothesis and \eqref{eq5} we deduce that $U_{k(n+1)+1}\equiv Q^kU_{k+1}Q^{\frac{k(n-2)}{2}}\equiv U_{k+1} Q^{\frac{kn}{2}}\pmod{U_k}$. This completes the proof.
\end{proof}
\end{lemma}
For the sake of brevity, it will be convenient to present the previous results without taking into consideration the parity of $n$. Therefore, we introduce a function $a(x)$ defined as \footnote{One may take for example $a(x)=\cos^2(\frac{\pi x}{2})$.}
\begin{gather*}a(x)=\begin{cases} 
   0, & \text{if}\quad x\equiv 1\pmod{2}, \\
      1, & \text{if}\quad x\equiv 0\pmod{2}.\\
   \end{cases}\end{gather*}
 In this way, we obtain the following Corollary.
\begin{corollary}\label{cor}
The following hold true:
\begin{equation}\label{eq6}\frac{U_{kn}}{U_k}\equiv n\left(U_{k+1}\right)^{a(n)}\cdot Q^{\frac{k(n-a(n)-1)}{2}}\pmod{U_k}.\end{equation}
\begin{equation}\label{eq7}U_{kn+1}\equiv \left(U_{k+1}\right)^{1-a(n)}\cdot Q^{\frac{k(n+a(n)-1)}{2}}\pmod{U_k}.\end{equation}
\end{corollary}
The previous Corollary is interesting on its own; we briefly describe two applications:\\
 From \eqref{eq6} we easily obtain that
\begin{gather*}U_{kn}=U_k\cdot n\left(U_{k+1}\right)^{a(n)}\cdot Q^{\frac{k(n-a(n)-1)}{2}}\pmod{U_k^2}.\end{gather*}
From this, it is straightforward to see that if $n^i\mid U_k$, then $n^{i+1}\mid U_{kn}$. This result has already appeared several times in the literature (for example, in \S 3 of \cite{2}, \emph{``The law of repetition''}), especially when $n$ is a prime number.\\
In addition, from the identity $U_{m+n}=U_{m+1}U_n-QU_mU_{n-1}$ for $m=kn$ and $n=r$ , we obtain that $U_{kn+r}=U_{kn+1}U_{r}-QU_{kn}U_{r-1}$. From \eqref{eq7} and the fact that $U_k\mid U_{kn}$, we deduce the more general congruence: 
\begin{gather*}U_{kn+r}\equiv U_r \left(U_{k+1}\right)^{1-a(n)}\cdot Q^{\frac{k(n+a(n)-1)}{2}}\pmod{U_k}.\end{gather*}
\section{Main results.}
In this section, we state our main results and provide some key proof techniques and insights. We recall that the Lucas sequence of the first kind is a strong divisibility sequence, which means that $\gcd(U_n, U_k)=U_{\gcd(n, k)}$.
\begin{theorem}\label{theor}Let $k\in\mathbb{N}^*$ and $\gcd(n, k)=d$. Then
\begin{gather*}\frac{U_{kn}}{U_k}\equiv\begin{cases} 
      dU_{d+1}\cdot \frac{U_n}{U_d}\cdot Q^{\frac{kn-k-n}{2}}\pmod{U_n}, &\text{if}\quad n, k\quad \text{are even and}\quad \frac{n}{d}, \frac{k}{d} \quad \text{are odd}, \\
     d\cdot \frac{U_n}{U_d}\cdot Q^{\frac{kn-k-n+d}{2}}\pmod{U_n}, & \text{otherwise}.\\
   \end{cases}\end{gather*} 
\begin{proof}
Suppose that 
\begin{equation}\label{eq8}\frac{U_{kn}}{U_k}=r+M\cdot U_n\end{equation}
where $0\leq r\leq U_n-1$ and $M\in\mathbb{N}$. This is equivalent to $U_{kn}=rU_k+MU_nU_k$. Since $U_n\mid U_{kn}$ and $U_n\mid kMU_nU_k$, we obtain $U_n\mid rU_k$. But $\gcd(U_n, U_k)=U_d$, which implies that $\gcd(\frac{U_n}{U_d}, \frac{U_k}{U_d})=1$. Thus, $\frac{U_n}{U_d}\mid r$, so $r=c\cdot \frac{U_n}{U_d}$ for some $c\in\mathbb{N}$. Hence, 
\begin{equation}\label{eq9}\frac{U_{kn}}{U_k}=c\cdot \frac{U_n}{U_d}+MU_n.\end{equation}  
From \eqref{eq7} for $k=d$ and $n=\frac{k}{d}$, we obtain 
\begin{equation}\label{eq10}U_{k+1}\equiv(U_{d+1})^{1-a(k/d)}\cdot Q^{\frac{d(k/d+a(k/d)-1)}{2}}\pmod{U_k}.\end{equation}
Thus, \eqref{eq6} implies that 
\begin{equation}\label{eq11}\frac{U_{kn}}{U_k}\equiv n(U_{d+1})^{a(n)(1-a(k/d))}\cdot Q^{\frac{da(n)a(k/d)-da(n)+kn-k}{2}}\pmod{U_k}.\end{equation}
Also from \eqref{eq6}, for $k=d$ and $n=\frac{n}{d}$, we obtain 
\begin{equation}\label{eq12}\frac{U_n}{U_d}\equiv \frac{n}{d}(U_{d+1})^{a(n/d)}\cdot Q^{\frac{d(n/d-a(n/d)-1)}{2}}\pmod{U_d}.\end{equation}
Since $n\equiv k\equiv 0\pmod{d}$, then $U_n\equiv U_k\equiv 0\pmod{U_d}$. Thus, \eqref{eq9}, \eqref{eq11} and \eqref{eq12} imply that
\begin{gather*}n(U_{d+1})^{a(n)(1-a(k/d))}\cdot Q^{\frac{da(n)a(k/d)-da(n)+kn-k}{2}}\equiv c\frac{n}{d}(U_{d+1})^{a(n/d)}\cdot Q^{\frac{d(n/d-a(n/d)-1)}{2}}\pmod{U_d} .\end{gather*}
Furthermore, we observe \footnote{ From definition, $U_{d+1}=P\cdot U_{d}-Q\cdot U_{d-1}$. If $\delta\mid U_d$ and $\delta\mid Q$, then $\delta\mid U_{d+1}$. Since  $\gcd(d, d+1)=1$, then $\gcd(U_d, U_{d+1})=1$, which implies that $\delta=1$.} that $\gcd(Q, U_d)=1$ and $\gcd(U_d, U_{d+1})=1$. As a consequence, we can rewrite the previous congruence in the form:
\begin{equation}\label{eq13}n(U_{d+1})^{a(n)(1-a(k/d))-a(n/d)}\cdot Q^{\frac{da(n)a(k/d)-da(n)+da(n/d)+kn-k-n+d}{2}}\equiv c\frac{n}{d}\pmod{U_d}.\end{equation}
Moving forward, we multiply both sides of \eqref{eq9} by $\frac{U_k}{U_n}$. This yields
\begin{gather*}\frac{U_{kn}}{U_n}=c\frac{U_k}{U_d}+MU_k.\end{gather*}
Following the same steps, (by exchanging the roles of $k$ and $n$) we conclude that 
\begin{equation}\label{eq14}\frac{U_{kn}}{U_n}\equiv k(U_{d+1})^{a(k)(1-a(n/d))}\cdot Q^{\frac{da(k)a(n/d)-da(k)+kn-n}{2}}\pmod{U_d}.\end{equation}
and 
\begin{equation}\label{eq15}\frac{U_k}{U_d}\equiv \frac{k}{d}(U_{d+1})^{a(k/d)}\cdot Q^{\frac{d(k/d-a(k/d)-1)}{2}}\pmod{U_d}.\end{equation}
Again, since $\gcd(Q, U_d)=\gcd(U_d, U_{d+1})=1$, we deduce that 
\begin{equation}\label{eq16}k(U_{d+1})^{a(k)(1-a(n/d))-a(k/d)}\cdot Q^{\frac{da(k)a(n/d)-da(k)+da(k/d)+kn-k-n+d}{2}}\equiv c\frac{k}{d}\pmod{U_d}.\end{equation}
The previous results indicate that the rest of the proof falls naturally into four cases: 
\begin{case1*}\emph{The number $n$ is odd.}\\
This implies that $d$ is odd (and so is $\frac{n}{d}$). Hence, $k$ and $\frac{k}{d}$ have the same parity. Therefore,
 \begin{flalign*}&a(n)(1-a(k/d))-a(n/d)=0,\\
 &da(n)a(k/d)-da(n)+da(n/d)=0,\\
 &a(k)(1-a(n/d))-a(k/d)=a(k)-a(k/d)=0,\\
 &da(k)a(n/d)-da(k)+da(k/d)=d(-a(k)+a(k/d))=0.\end{flalign*}
\end{case1*}
\begin{case2.a*}\emph{The number $n$ is even and $k$ is odd.}\\
This implies that $d$ and $\frac{k}{d}$ are odd but $\frac{n}{d}$ is even. Therefore, 
\begin{flalign*}&a(n)(1-a(k/d))-a(n/d)=a(n)-a(n/d)=0,\\
 &da(n)a(k/d)-da(n)+da(n/d)=-d(a(n)-a(n/d))=0,\\
&a(k)(1-a(n/d))-a(k/d)=0,\\
&da(k)a(n/d)-da(k)+da(k/d)=0.\end{flalign*}
\end{case2.a*}
\begin{case2.b*} \emph{Both $n$ and $k$ are even and $\frac{n}{d}, \frac{k}{d}$ have different parity.}\\
Since $\frac{n}{d}, \frac{k}{d}$ have different parity, then $a(n/d)+a(k/d)=1$. Therefore,
\begin{flalign*}&a(n)(1-a(k/d))-a(n/d)=1-\left(a(k/d)+a(n/d)\right)=0,\\
&da(n)a(k/d)-da(n)+da(n/d)=d(a(k/d)+a(n/d)-1)=0,\\
&a(k)(1-a(n/d))-a(k/d)=1-\left(a(n/d)+a(k/d)\right)=0,\\
&da(k)a(n/d)-da(k)+da(k/d)=d(a(n/d)-1+a(k/d))=0.\end{flalign*}\end{case2.b*} 
In all these cases, we get 
\begin{equation}\label{eq17}k\cdot Q^{\frac{kn-k-n+d}{2}}\equiv c\frac{k}{d}\pmod{U_d}.\end{equation}
and
\begin{equation}\label{eq18}n\cdot Q^{\frac{kn-k-n+d}{2}}\equiv c\frac{n}{d}\pmod{U_d}.\end{equation}
For the convenience of notation, we let $Q^{\frac{kn-k-n+d}{2}}=x$. Congruences \eqref{eq17} and \eqref{eq18} imply that $\frac{k}{d}(dx-c)\equiv \frac{n}{d}(dx-c)\equiv 0\pmod{U_d}$. We claim that $dx-c\equiv 0\pmod{U_d}$. Indeed, since $\gcd(\frac{n}{d}, \frac{k}{d})=1$, Bezout's identity implies that there are $s, q \in \mathbb{Z}$ such that $s\cdot \frac{n}{d}+q\cdot \frac{k}{d}=1$. Therefore,
\begin{gather*}s\cdot\frac{n}{d}(dx-c)+q\cdot\frac{k}{d}(dx-c)=dx-c.\end{gather*}
Thus, $dx-c\equiv 0\pmod{U_d}$. Consequently, $U_n\frac{dx-c}{U_d}$ is an integer multiple of $U_n$ which is equivalent to $c\frac{U_n}{U_d}\equiv dx\frac{U_n}{U_d}\pmod{U_n}$. Since $\frac{U_{kn}}{U_k}\equiv c\frac{U_n}{U_d}\pmod{U_n}$, it follows that
\begin{gather*}\frac{U_{kn}}{U_k}\equiv d\cdot \frac{U_n}{U_d}\cdot Q^{\frac{kn-k-n+d}{2}}\pmod{U_n}.\end{gather*}
\begin{case2.c*}\emph{Both $n$ and $k$ are even and $\frac{n}{d}, \frac{k}{d}$ are odd\footnote{There is no other case to consider, since $\frac{n}{d}, \frac{k}{d}$ cannot be both even since they are coprime.}.}\\
In this case, we observe that 
\begin{flalign*}&a(n)(1-a(k/d))-a(n/d)=1,\\
&da(n)a(k/d)-da(n)+da(n/d)=da(k/d)-d+da(n/d)=-d,\\
&a(k)(1-a(n/d))-a(k/d)=1,\\
&da(k)a(n/d)-da(k)+da(k/d)=-d.\end{flalign*}\end{case2.c*} Hence, from \eqref{eq13} and \eqref{eq16} we deduce that 
\begin{gather*}kU_{d+1}\cdot Q^{\frac{kn-k-n}{2}}\equiv c\frac{k}{d}\pmod{U_d}.\end{gather*} and
\begin{gather*}nU_{d+1}\cdot Q^{\frac{kn-k-n}{2}}\equiv c\frac{n}{d}\pmod{U_d}.\end{gather*}
Let $U_{d+1}\cdot Q^{\frac{kn-k-n}{2}}=y$. Using the same technique as before, we obtain that $U_n\frac{dy-c}{U_d}$ is an integer multiple of $U_n$, which implies $\frac{U_{kn}}{U_k}\equiv c\frac{U_n}{U_d}\equiv dy\frac{U_n}{U_d}\pmod{U_n}$. It follows that
\begin{gather*}\frac{U_{kn}}{U_k}\equiv dU_{d+1}\cdot \frac{U_n}{U_d}\cdot Q^{\frac{kn-k-n}{2}}\pmod{U_n}.\end{gather*}
This completes the proof. 
\end{proof}
\end{theorem}
We observe that if $Q=\pm 1$, we can simplify the shape of the congruences which appear in Corollary \eqref{cor} and Theorem\eqref{theor}. There is a plethora of examples of Lucas sequences $U_n(P, \pm 1)$. In order to limit the length of this paper, we chose the Fibonacci sequence $(F_n)_n$ which is obtained from $U_n(1, -1)$ and the Mersenne numbers, which are numbers of the form $U_n(3, 2)=2^n-1$ to demonstrate the application of our results. 
\begin{theorem} The following hold true:
 \begin{flalign}&\label{eq19}F_{kn+1}\equiv\begin{cases} 
      F_{k+1}(-1)^{\frac{k(n-1)}{2}}\pmod{F_k}, &\text{if}\quad n\equiv 1\pmod{2}, \\
      (-1)^{\frac{kn}{2}}\pmod{F_k}, &\text{if}\quad n\equiv 0\pmod{2}.\\
   \end{cases}\\
&\label{eq20}\frac{F_{kn}}{F_k}\equiv\begin{cases} 
      n\cdot (-1)^{\frac{k(n-1)}{2}}\pmod{F_k}, &\text{if}\quad n\equiv 1\pmod{2}, \\
      nF_{k+1}\cdot (-1)^{\frac{k(n-2)}{2}}\pmod{F_k}, &\text{if}\quad n\equiv 0\pmod{2}.\\
\end{cases}\\
&\label{eq21}\frac{F_{kn}}{F_k}\equiv\begin{cases}
d\cdot \frac{F_n}{F_d}\cdot F_{d+1} (-1)^{\frac{kn-k-n}{2}}\pmod{F_n}, &\text{if}\quad n, k\quad \text{are even and}\quad \frac{n}{d}, \frac{k}{d} \quad \text{are odd},\\
d\cdot \frac{F_n}{F_d}\cdot (-1)^{\frac{kn-k-n+d}{2}}\pmod{F_n}, & \text{otherwise}.
\end{cases}\\
&\label{eq22}\frac{2^{kn}-1}{2^k-1}\equiv d\cdot \frac{2^n-1}{2^d-1}\pmod{2^n-1}.
\end{flalign}\\
\begin{proof}
The first three parts follow immediately from Corollary \eqref{cor} and Theorem\eqref{theor} if we substitute $Q=-1$. In order to prove \eqref{eq22}, we substitute $Q=2$ into Theorem \eqref{theor} and treat separately the two cases:\\
\emph{If $n, k$ are even and $\frac{n}{d}, \frac{k}{d}$ are odd}:\\
It is evident that $d$ and $\frac{n}{d}+\frac{k}{d}$ are even. Thus,  
\begin{gather*}\frac{kn-k-n}{2}=d\cdot\frac{(d\cdot \frac{n}{d}\frac{k}{d}-(\frac{k}{d}+\frac{n}{d}))}{2}\equiv 0\pmod{d}.\end{gather*}
This implies that $2^{\frac{kn-k-n}{2}}\equiv 1\pmod{2^d-1}$. Furthermore, $U_{d+1}=2^{d+1}-1\equiv 1 \pmod{2^d-1}$. We conclude that
\begin{gather*}\frac{2^{kn}-1}{2^k-1}\equiv d(2^{d+1}-1)\frac{2^n-1}{2^d-1}2^{\frac{kn-k-n}{2}}\equiv d\cdot \frac{2^n-1}{2^d-1}\pmod{2^n-1}.\end{gather*}
\emph{Otherwise}:\\
Observe that 
\begin{gather*}\frac{kn-k-n+d}{2}=d\cdot \frac{d\frac{n}{d}\frac{k}{d}-\frac{k}{d}-\frac{n}{d}+1}{2}.\end{gather*}
If $d$ is odd, then obviously $d\cdot \frac{d\frac{n}{d}\frac{k}{d}-\frac{k}{d}-\frac{n}{d}+1}{2}\equiv 0\pmod{d}$. \\
If $d$ is even, then both $n$ and $k$ must be even and $\frac{k}{d}, \frac{n}{d}$ must have different parity. (We already considered above the case where both are odd.) Therefore, $d\frac{n}{d}\frac{k}{d}-\frac{k}{d}-\frac{n}{d}+1\equiv 0\pmod{2}$, which implies that $d\cdot \frac{d\frac{n}{d}\frac{k}{d}-\frac{k}{d}-\frac{n}{d}+1}{2}\equiv 0\pmod{d}$. In any case, $2^{\frac{kn-k-n+d}{2}}\equiv 1\pmod {2^d-1}$ holds true. It follows that 
\begin{gather*}\frac{2^{kn}-1}{2^k-1}\equiv d\cdot \frac{2^n-1}{2^d-1}2^{\frac{kn-k-n+d}{2}}\equiv d\cdot \frac{2^n-1}{2^d-1}\pmod{2^n-1}.\end{gather*}
This completes the proof.
\end{proof}
\end{theorem}
We would like to note that the behavior of $\frac{F_{kn}}{F_k}$ and $F_{kn+1}$ modulo $F_k$ goes back at least as far as \cite{4}, where the author obtains that $\frac{F_{kn}}{F_k}\equiv n(F_{k+1})^{n-1}\pmod{F_k}$ and $F_{kn+1}\equiv (F_{k+1})^n\pmod{F_k}$. Since $n(F_{k+1})^{n-1}$ and $(F_{k+1})^n$ are much larger than the moduli we obtained, we consider \eqref{eq19} and \eqref{eq20} to be much more illuminating.
\subsection{Criteria of Primality.}
\begin{theorem}\label{theorem.}The congruence \begin{gather*}\sum_{k=1}^{n-1}\frac{2^{kn}-1}{2^k-1}\equiv 0\pmod{2^n-1}\end{gather*} holds true, if and only if $n$ is prime.
\begin{proof}

From \eqref{eq22}, we observe that 
\begin{gather*}\sum_{k=1}^{n-1}\frac{2^{kn}-1}{2^k-1}\equiv \sum_{\substack{\gcd(n, k)=d\\ 1\leq k<n}}d\cdot \frac{2^n-1}{2^d-1}\pmod{2^n-1}.\end{gather*}
The number of positive integers $k$ less than $n$ with $\gcd(n, k)=d$, is exactly $\varphi(\frac{n}{d})$. Thus,
\begin{equation}\label{eq23}\sum_{k=1}^{n-1}\frac{2^{kn}-1}{2^k-1}\equiv (2^n-1)\sum_{\substack{d\mid n\\ 1\leq d<n}}\frac{d\cdot\varphi(\frac{n}{d}) }{2^d-1}\pmod{2^n-1}.\end{equation}
It suffices to prove that the last sum is an integer if and only if $n$ is prime.\\
If $n$ is prime, there is only one divisor $d<n$ of $n$, namely $d=1$ and the sum is equal to $\frac{1\cdot \varphi(n)}{2^1-1}=n-1$.\\
If $n$ is composite, we denote by $D$ the largest proper divisor of $n$. From Theorem \eqref{theore}, we know that $2^D-1$ has a primitive prime divisor $p$, unless $D=6$. It is relatively easy to show that $n=12$ is the only number with $D=6$, and $\sum_{k=1}^{11}\frac{2^{k\cdot 12}-1}{2^k-1}\equiv 3354\not\equiv0\pmod{2^{12}-1}$. Therefore, we suppose that $D\neq 6$.  We denote by $L$ the least common multiple of all the elements of $\{2^d-1, 1\leq d<n\}$ and $d_1, d_2, \ldots, D$ all the proper divisors of $n$. Hence, 
\begin{equation}\label{eq24}\sum_{\substack{d\mid n\\ 1\leq d<n}}\frac{d\cdot\varphi(\frac{n}{d}) }{2^d-1}=\frac{d_1\varphi(\frac{n}{d_1})\left(\frac{L}{2^{d_1}-1}\right)+d_2\varphi(\frac{n}{d_2})\left(\frac{L}{2^{d_2}-1}\right)+\ldots +D\varphi(\frac{n}{D})\left(\frac{L}{2^{D}-1}\right)}{L}.\end{equation}
 
We briefly show that $p\nmid D\varphi(\frac{n}{D})$. From Fermat's little theorem we know that $p\mid 2^{p-1}-1$, which implies that $p>D$. Hence, $p\nmid D$. Furthermore, since $D$ is the largest proper divisor of $n$, then $\frac{n}{D}$ is the least prime divisor of $n$. This implies that $\varphi(\frac{n}{D})= \frac{n}{D}-1$ and $\frac{n}{D}\leq D$. Thus, $\varphi(\frac{n}{D})<D<p$, which shows that $p\nmid \varphi(\frac{n}{D})$. 
By this reasoning, the sum in \eqref{eq24} is not an integer, since the denominator and every summand in the numerator is divisible by $p$, except $D\varphi(\frac{n}{D})\left(\frac{L}{2^{D}-1}\right)$. This completes the proof.
\end{proof}
\end{theorem}
 Before proceeding any further, we would like to mention two well-known theorems which will be needed in the proof of our last main result.
In 1913, Carmichael \cite{5} showed that almost all terms of a Lucas sequence of the first kind have a \emph{primitive prime divisor}. In particular he proved the following Theorem:
\begin{theorem} \label{theore}Suppose that $P^2-4Q>0$. If $n\neq1, 2, 6$, $U_n$ has a prime divisor $p$, such that $p\nmid U_m$, for every $m<n$. The only exception is $U_{12}(1, -1)=F_{12}=144$.
\end{theorem} 
We will also make use of the following theorem which can be found in various contexts, such as \cite{2}:
\begin{theorem}\label{theorem}
If $p$ is an odd prime, $p\neq 5$, then $p\mid F_{p-\legendre{5}{p}}$, where $\legendre{a}{p}$ denotes the Legendre symbol.
\end{theorem}
Roughly speaking, this theorem shows that if $p$ is the primitive prime divisor of $F_n$, then $p\geq n-1$. To see why this is true, observe that $\legendre{5}{p}=\pm 1$ if $p\neq 5$, and $\legendre{5}{p}=0$ if $p=5$.

\begin{theorem}Suppose that $n\equiv 1\pmod{4}$ and $n\not\in \{9, 25\}$. Then
\begin{equation}\label{eq25}\sum_{k=1}^{n-1}\frac{F_{kn}}{F_k}\equiv 0\pmod{F_n}\end{equation}
holds true, if and only if $n$ is prime.
\begin{proof}
This theorem may be proved in much the same way as Theorem\eqref{theorem.}. From \eqref{eq21}, we get
\begin{gather*}\sum_{k=1}^{n-1}\frac{F_{kn}}{F_k}\equiv F_n\cdot \sum_{\substack{\gcd(n, k)=d\\ 1\leq d<n}}\frac{d}{F_d}(-1)^{\frac{kn-k-n+d}{2}}\pmod{F_n}.\end{gather*}
It will be convenient to substitute $-1$ with $i^2$, the imaginary unit. In this way, we may write $(-1)^{\frac{kn-k-n+d}{2}}=i^{d-1}(i^{n-1})^{(k-1)}$. Consequently, it suffices to prove that
\begin{gather*}\sum_{\substack{\gcd(n, k)=d\\ 1\leq d<n}}\frac{d}{F_d}i^{d-1}(i^{n-1})^{(k-1)}=\sum_{\substack{d\mid n\\ 1\leq d<n}}\frac{di^{d-1}}{F_d}\sum_{\substack{\gcd(n, k)=d}}(i^{n-1})^{(k-1)}\end{gather*}
is an integer if and only if $n$ is prime.
If $n\equiv 1\pmod{4}$, then $i^{n-1}=1$, which implies that $\sum_{\substack{\gcd(n, k)=d}}(i^{n-1})^{(k-1)}$ is equal to $\varphi(\frac{n}{d})$. It is evident that if $n$ is prime, the sum is an integer, namely $\varphi(n)=n-1$. Suppose that $n$ is composite, $D$ is its largest proper divisor and $p$ is the primitive prime divisor of $F_D$ \footnote{Since $n$ is odd, we do not need to consider the exceptions mentioned in Theorem \eqref{theore}.}. As is clear from the previous arguments, it will be enough to show that 
\begin{gather*}p\nmid Di^{D-1}\varphi\left(\frac{n}{D}\right).\end{gather*}
Recall that from Theorem \eqref{theorem}, we know that $p\geq D-1$. Thus, if $p\mid D$, then $p=D$. This implies that $p\mid F_p$, which is possible only if $p=D=5$. The only composite number congruent to $1$ modulo $4$, having $D=5$, is $25$. From hypothesis $n\neq 25$, thus, $p\nmid D$.\\ 
Suppose that $p\mid\varphi(\frac{n}{D})$. Then, 
\begin{gather*}p\leq \varphi\left(\frac{n}{D}\right)=\frac{n}{D}-1\leq D-1\leq p.\end{gather*}
This shows that $p=D-1$ and $\frac{n}{D}-1=D-1$. The last equality implies that $n=D^2$, which is possible only if $D$ is prime. Hence, $p, D$ are consecutive primes which yields $p=2, D=3, n=9$. From hypothesis $n\neq 9$, thus, $p\nmid  \varphi(\frac{n}{D})$.
This completes the proof.
\end{proof}
\end{theorem}
\begin{remark}
We pause now to record an observation that follows readily from the computations made above; if $n\equiv 3\pmod{4}$, then \eqref{eq25} holds true, regardless of whether $n$ is prime. To see this, observe that $i^{n-1}=-1$, thus
\begin{gather*}\sum_{\substack{\gcd(n, k)=d}}(i^{n-1})^{(k-1)}= -\sum_{\gcd(n, k)=d}(-1)^{d\cdot \frac{k}{d}}=-\sum_{\gcd(\frac{n}{d}, \frac{k}{d})=1}(-1)^{\frac{k}{d}}.\end{gather*} The last equality holds, since $d$ is odd. In addition, we note that $\frac{k}{d}$ denotes a number which is coprime to $\frac{n}{d}$. If $\gcd(\frac{k}{d}, \frac{n}{d})=1$, then $\gcd(\frac{n}{d}-\frac{k}{d}, \frac{n}{d})=1$. But $\frac{n}{d}-\frac{k}{d}$ and $\frac{k}{d}$ have opposite parity since their sum is $\frac{n}{d}\equiv 1 \pmod{2}$, which means that $(-1)^{\frac{n}{d}-\frac{k}{d}}+(-1)^{\frac{k}{d}}=0$. Hence, for every $d$ we have a sum of pairs which vanishes, and the result follows.
\end{remark}
\makeatletter
\renewcommand{\@biblabel}[1]{[#1]\hfill} 
\makeatother

\end{document}